\documentclass[reqno]{amsart}
\usepackage{amssymb}

\def\S{\mathbb S}
\def\N{\mathbb N}
\newcommand{\R}{\mathbb R}

\newcommand\rr[1]{\mbox{(\ref{#1})}}

\newcommand{\sect}[1]{\section{#1}\setcounter{equation}{0}}

\newtheorem{theorem}{Theorem}[section]

\newtheorem{proposition}[theorem]{Proposition}

\newtheorem{remark}[theorem]{Remark}
\newtheorem{example}[theorem]{Example}

\def\qed{\hfill$\blacksquare$\vskip .2cm}

\begin{document}
\author{Fred W. Huffer, Jayaram Sethuraman, and Sunder Sethuraman}

\address{\noindent Department of Statistics,
Florida State University, Tallahassee, FL \ 32306.
\newline
e-mail:  \rm \texttt{huffer@stat.fsu.edu}
}

\address{\noindent Department of Statistics,
Florida State University, Tallahassee, FL \ 32306.
\newline
e-mail:  \rm \texttt{sethu@stat.fsu.edu}
}

\address{\noindent Department of Mathematics,
396 Carver Hall, Iowa State University, Ames, IA \ 50011, USA.
\newline
e-mail:  \rm \texttt{sethuram@iastate.edu}}

\title[Bernoulli strings and Poisson processes]
{A study of counts of Bernoulli strings via conditional Poisson
processes}

\begin{abstract}
  A sequence of random variables, each taking values $0$ or $1$,
  is called a Bernoulli sequence. We say that a string of length
  $d$ occurs, in a Bernoulli sequence, if a success is followed
  by exactly $(d-1)$ failures before the next success.  The
  counts of such $d$-strings are of interest, and
in specific independent Bernoulli
  sequences are known to correspond to asymptotic $d$-cycle
  counts in random permutations.

  In this note, we give a new framework, in terms of conditional
Poisson processes, which allows for a quick characterization of
  the joint distribution of the counts of all $d$-strings,
in a general class of Bernoulli sequences,
  as certain mixtures of the product of Poisson measures. In
  particular, this general class includes all Bernoulli sequences
  considered in the literature, as well as a host of new sequences.
\end{abstract}



\thanks{Research partially supported by ARO-W911NF-04-1-0333, NSA-H982300510041,
  and NSF-DMS-0504193.}

\maketitle

\sect{Introduction}
\label{intro_sect}
In this note, we study the joint distribution of
the counts of certain $d$-strings of all orders $d>1$ arising in
Bernoulli sequences.
Previous work has used several different methods, including
combinatorial, factorial moment, and
P\'olya and Hoppe urn model methods to identify
the joint count distribution with respect to
a class of independent Bernoulli sequences.
In this context, our main contribution is to introduce a new framework, using
conditional Poisson processes, which allows for a concise derivation
of
the joint count distribution as a
mixture of the product of Poisson measures with respect to all
Bernoulli sequences considered before, as well as many others in a
general class,
including some dependent Bernoulli sequences.





A Bernoulli sequence ${\bf Y} = \{Y_n\}_{n\geq 1}$ is a sequence
of $\{0,1\}$-valued random variables. For $d\geq 1$, we say that
a $d$-string occurs if a $1$ is followed by exactly $(d-1)$ $0$'s
before the next $1$ in the Bernoulli sequence.  Specifically, a
$d$-string occurs at time $n\geq 1$ if $Y_{n,d}=1$ where
$$Y_{n,d} \ {=}\  \left\{\begin{array}{rl}
Y_nY_{n+1}& \ {\rm for \ }d=1\\
Y_n(1-Y_{n+1})\cdots (1-Y_{n+d-1})Y_{n+d}&\ {\rm for \ }d\geq
2,\end{array}\right.
$$
that is, if $\langle Y_n,\dots,Y_{n+d}\rangle  =  \langle
1,\underbrace{0,\dots,0}_{d-1},1\rangle.$

Let $Z_d=\sum_{n\geq 1}Y_{n,d}$ be the count of all $d$-strings,
for $d\geq 1$, and ${\bf Z} = \langle Z_d: d\geq 1\rangle$ be
the ``count vector'' of strings.
[In general, ${\bf Z}$ may have
divergent components, but for the Bernoulli sequences
considered in this article it is easily shown (by taking expectations)
that all components
$Z_k$ are finite with probability $1$.]

In this notation,
the general problem is to understand the distribution of ${\bf Z}$ and its
connection to the underlying sequence ${\bf Y}$.
Aside from the problem's basic
interest, $d$-strings and their counts from specific independent
Bernoulli sequences have interpretations with
respect to random permutations, record values, Bayesian
nonparametrics, and species allocation models through Ewens
sampling formula.

We will use ``$\stackrel{d}{=}$'' to signify
``equals in distribution,'' and ${\mathcal L}(X)$ to denote the law or
distribution of the random variable $X$.  Denote also ${\rm Po}(\lambda)$ as
the Poisson measure on $\R$ with intensity $\lambda$, and $I(B)$
as the indicator of a set $B$.

\begin{example}
\label{ex_1}
\rm    Let $\S_n = \{1,2,\ldots,n\}$, and consider the Feller algorithm to generate a permutation
$\pi:\S_n\rightarrow \S_n$
uniformly among the $n!$ choices (cf. Feller~(1945)):
\begin{itemize}
\item[1.] Draw an element uniformly from $\S_n$, and call it
  $\pi(1)$.  If $\pi(1)=1$, a $1$-cycle is completed.  If
  $\pi(1)\neq 1$, make another draw uniformly from
  $\S_n\setminus\{\pi(1)\}$, and call it $\pi(\pi(1))$.
  Continue drawing from $\S_n\setminus\{\pi(1),
  \pi(\pi(1))\},\ldots$
naming them $\pi(\pi(\pi(1)))$, and so on,
until a cycle (of some length) is
  finished.

\item[2.] From the elements left in $\S_n\setminus\{\pi(1),\pi(\pi(1)), \ldots,1\}$
after the first cycle is
completed, follow the process in step $1$ with the smallest remaining
number taking the role of ``$1$'' to finish a second cycle.  Repeat
until all elements of $\S_n$ are exhausted.
\end{itemize}

Let $I^{(n)}_k$ be the indicator that a cycle is
completed at the $k$th Feller draw from $\S_n$.  A moment's thought
convinces
that $\{I^{(n)}_k\}_{k=1}^n$ are independent Bernoulli random
variables with $P(I^{(n)}_k=1)=1/(n-k+1)$ as, independent of the past,
exactly one choice at time $1\leq k\leq n$ from the remaining
$n-k+1$ members left in $\S_n$ completes the cycle.  Denote
$C^{(n)}_k$ as the number of $k$-cycles in $\pi$,
$$C^{(n)}_k  = \left\{\begin{array}{rl}
I^{(n)}_1 + \sum_{i=1}^{n-1} I^{(n)}_iI^{(n)}_{i+1}&  {\rm for \ }k=1\\
 \prod_{l=1}^{k-1}(1-I^{(n)}_l)I^{(n)}_k +
 \sum_{i=1}^{n-k}I^{(n)}_i\prod_{l=i+1}^{i+k-1}(1-I^{(n)}_l)I^{(n)}_{i+k}&  {\rm for \ }2\leq k\leq n.\end{array}\right.$$

Now let ${\bf Y}$ be the independent sequence where
$P(Y_k=1)=1/k$
for $k\geq 1$, so that $Y_k\stackrel{d}{=}I^{(n)}_{n-k+1}$ for $1\leq
k\leq n$.  Then, as $Y_n$, and $Y_{n-k+1}\prod_{l=n-k+2}^n(1-Y_l)$ for
$2\leq k\leq n$ all vanish in probability as $n\uparrow \infty$, we
conclude for each $k\geq 1$ that
$\lim_{n\rightarrow \infty}C^{(n)}_k  \stackrel{d}{=}  Z_k$.

Finally, as is well-known, the asymptotic cycle counts $\{\lim_n
C^{(n)}_k\}_{k\geq 1}$ are distributed as independent Poisson
random variables with respective means $1/k$ for $k\geq 1$ (cf.
Kolchin~(1971)).  Hence, ${\bf Z}\stackrel{d}{=}\prod_{k\geq
  1}{\rm Po}(1/k)$. [Example~\ref{ex_3}, in section 2,
gives a derivation in our Poisson process framework.  See also
Arratia-Barbour-Tavar\'e~(1992,\ 2003) for more discussion with
Ewens sampling formula.]
\end{example}
\begin{example}
  \label{ex_2}
  \rm Consider the standard nonparametric problem of estimating
  the unknown distribution function $F$ from independent and
  identically distributed observations $\{X_i\}_{i\geq 1}$. A
  Bayesian may place on $F$
a Dirichlet prior with parameters $a\mu$ where $a>0$ and $\mu$ is
  a non-atomic probability measure.

Let $Y_1=1$ and for $n\ge
  2$ define $Y_n =1$ if $X_n$ is a new observation, that is if $X_n
  \not \in \{X_1,\dots,X_{n-1}\}$, and $Y_n=0$ otherwise. Then, it can be shown that ${\bf
  Y}$ is an independent Bernoulli sequence
with $P(Y_n=1) = a/(a+n-1)$ for $n\geq 1$ and that 
$(\log n)^{-1}\sum_{i=1}^n Y_i  \rightarrow  a$ a.s. The latter
result can be interpreted in terms of counts of strings in this
Bernoulli sequence.
See Korwar-Hollander~(1973) for
more details, and also Ghosh-Ramamoorthi (2003).

\end{example}

In the literature, to our knowledge, only
the count vectors of the following class of underlying
independent Bernoulli sequences have been investigated. Denote
the independent Bernoulli sequence ${\bf Y}$ where
$P(Y_n=1)=a/(a+b+n-1)$ for $n\geq 1$ as ${\bf Y}={\rm
  Bern}(a,b)$.  The case $a=1$, $b=0$ is Example~\ref{ex_1} (see
also Arratia-Tavar\'e~(1992)). The case $a>0$, $b=0$ is
Example~\ref{ex_2}. For this case, Arratia-Barbour-Tavar\'e~(1992)
observe that the associated ${\bf Z}\stackrel{d}{=}\prod_{k\geq
  1}{\rm Po}(a/k)$ through connections with Ewens sampling
formula.  When $a=1$, $b>0$, Sethuraman-Sethuraman~(2004),
employing factorial moments, show that, given the value $x_0$ of
a Beta$(b,1)$ random variable, ${\bf
  Z}\stackrel{d}{=}\prod_{k\geq 1} {\rm Po}((1-x_0^k)/k)$.  Such a
distribution will be called a ``mixture of independent Poisson
factors.'' When $a>0$ and $b> 0$, Holst (2007) extends further,
using P\'olya and Hoppe urns, and establishes that, given the
value $x_0$ of a Beta$(b,a)$ random variable, ${\bf
  Z}\stackrel{d}{=}\prod_{k\geq 1} {\rm Po}(a(1-x_0^k)/k)$,
again a mixture of independent Poisson factors.  We note also
that several interesting studies of $1$-strings preceded
some of the above work, e.g. an unpublished manuscript of
Diaconis, Chern-Hwang-Yeh~(2000), M\'ori~(2001),
Joffe-Marchand-Perron-Popadiuk~(2004), and references therein in these
and the above papers.

With this background,
our main
idea is that it is easier to study ${\bf Z}$
starting from an extrinsic
``conditional marked Poisson process model'' (CMPP) rather
 than directly from the Bernoulli sequence.  Namely, we prove
that
when the underlying
Bernoulli sequence ${\bf Y}$ is generated through a CMPP model,
the count vector ${\bf Z}$
is distributed as a mixture of independent Poisson factors in
terms of model parameters (Theorem~\ref{cmpp_thm}).
As remarked earlier,
the Poisson process techniques used
here are different from previous methods and allow quick derivations.
Perhaps interestingly, the sequences ${\bf Y}$ found in our model include many
dependent Bernoulli sequences (some explicit examples are in
section~\ref{dep_sect}).  However, the most general sequence
studied till now, the independent sequence ${\rm Bern}(a,b)$ with
$a>0$ and $b\geq 0$, can also be realized in our framework
(Proposition \ref{holst_prop}), yielding a new proof of its count
vector distribution.

Our conditional marked Poisson process model also yields a new
class of independent Bernoulli sequences which we call ${\rm
  Bern}_1(a,b)$.  Denote the independent Bernoulli sequence ${\bf
  Y}$ where $P(Y_1=1)=1$, and $P(Y_n=1)=a/(a+b+n-2)$ for $n\geq 2$ as
${\bf Y}={\rm Bern}_1(a,b)$.  The ${\rm Bern}_1(a,b)$ sequence
appends a $1$ to the ${\rm Bern}(a,b)$ sequence and picks up one
more $d$-string contributed by any leading $0$'s in ${\rm
  Bern}(a,b)$.
We show
that the distribution of the count vector ${\bf Z}$ for ${\rm Bern}_1(a,b)$
for $a>0,b \geq 1$ is a mixture of independent Poisson factors
(Proposition~\ref{new_prop}).  This result fails for $0 \leq b<1$, and in this case
even the
distribution of $Z_1$, the count of $1$-strings in ${\rm Bern}_1(a,b)$,
is not a mixture of Poisson distributions
(Proposition~\ref{cont_prop}). However, the distribution of ${\bf Z}$
in ${\rm Bern}_1(a,b)$ can be expressed through a recurrence relation
for all values of $b$ including $0\le b <1$
(Proposition~\ref{decomp_prop}).

 The plan of the article is to discuss
the CMPP model, and prove the main theorem in section ~\ref{cmpp_sect}.
In sections~\ref{holst_sect} and \ref{new_sect}, the main theorem
 is applied to independent sequences ${\rm
   Bern}(a,b)$ and ${\rm Bern}_1(a,b)$ respectively.  Last, in section
  \ref{dep_sect}, two explicit dependent Bernoulli sequences,
 arising from the CMPP model, are given.

\sect{CMPP models}
\label{cmpp_sect}
The following ``Poisson process''
derivation of the distribution of ${\bf Z}$ with
respect to ${\rm Bern}(1,0)$ (cf. Example \ref{ex_1}) motivates subsequent development.

\begin{example}
\label{ex_3}
\rm Consider the following
standard way to generate a ${\rm Bern}(1,0)$ sequence.
Let $\{\beta_i\}_{i\geq 1}$ be independent, identically distributed (iid)
Uniform$[0,1]$ random
variables, and define $Y_n = I(\beta_n {\rm\ is
  \ a \ record}),\ n\ge 1$.  R\`enyi's theorem shows that
$\{Y_n\}_{n\geq 1}$ are independent and $P(Y_n =1)=1/n$ for $n\geq 1$,
that is ${\bf Y} = {\rm Bern}(1,0)$.  Let $\{X_i\}_{i\geq 1}$ be the
record values among $\{\beta_i\}_{i \ge 1}$. Notice that the point
process $N$ on $[0,1]$ defined by $N(A) = \sum_{i\geq
  1}\delta_{X_i}(A)$ is a nonhomogeneous Poisson process on $[0,1]$
with intensity $1/(1-x)$ (cf.  Resnick~(1994)).  For each point $X_i$,
we can associate a Geometric$(1-X_i)$ variable $L_i$ (a ``mark'')
corresponding to the number of uniform random variables in
$\{\beta_i\}_{i\geq 1}$ to the next record.  Then, by thinning
decompositions, $Z_k= \sum_{i\geq 1}I(L_i=k) = \sum_{i\geq 1}
\delta_{X_i}([0,1])I(L_i=k)$ for $k\geq 1$ are independent Poisson
variables with respective means $\int_0^1 (1-x)^{-1}x^{k-1}(1-x)dx =
1/k$ for $k\geq 1$.
\end{example}

In a sense, the thrust of the following CMPP model and our main result
(Theorem \ref{cmpp_thm}) below
is to reverse the procedure in Example~\ref{ex_3}.
By beginning with a given Poisson process and
spacing variables, which themselves determine the count vector ${\bf
 Z}$, we then see what associated Bernoulli sequence ${\bf Y}$ arises.

Consider a sequence of random variables ${({\bf X,L})}
 =\{(X_i,L_i)\}_{i\geq 0}$ on $\R\times \N$ where $\N = \{1,2,\ldots\}$, and the
point process $N$ on $\R$ given by
$N(A) {=} \sum_{i\geq 1}\delta_{X_i}(A)$.  Let also $g:\R\rightarrow
 [0,\infty)$
 be a
probability density function (pdf), and for each $x\in \R$
$r(x,\cdot),q(x,\cdot):\N\rightarrow [0,1]$
be probability mass functions, and
$\lambda_{x}: \R \rightarrow [0,\infty)$ be an
intensity function.

Then, we say ${({\bf X,L})}$ is the conditional marked
Poisson process
${\mathcal M}(g,r,\lambda,q)$ if the following
hold:
\begin{itemize}
\item [{1.}] $X_0$ has pdf $g$,
\item [{ 2.}] conditional on $X_0=x_0$, $N$ is a nonhomogeneous Poisson
  process with intensity function $\lambda_{x_0}(\cdot)$,
\item [{ 3.}] $P(L_0=k|{\bf X}) = r(X_0,k)$ for $k\geq 1$, and
\item [{ 4.}] $P(L_n=k|{\bf X},L_0,L_1,\dots,L_{n-1}) = q(X_n,k)$
  for $k,n\geq 1$.
\end{itemize}

Let $L_0^*=L_0$, and $L^*_r=L^*_{r-1}+L_r$ for $r\geq 1$. We now
define a Bernoulli sequence ${\bf Y}$ based on $({\bf X,L})$ as
follows:  $Y_n=1$ if $n$ is of the form $L^*_r$ for some $r\geq
0$, and $Y_n=0$ otherwise.  Another way to say this is
\begin{equation}
\label{bernoulli_seq}
Y_n \ =\ \left\{\begin{array}{rl}
0 & {\rm \ when \ }  n< L_0^*, {\rm \ or \ }L_r^* < n < L_{r+1}^* {\rm \
  for \ }
r\geq 0\\
1 & {\rm \ when \ }  n=L_r^* {\rm \ for \ }r\geq
0.\end{array}\right.
\end{equation}
Then, the count vector ${\bf Z}$ is given by
\begin{equation}
\label{count_rep}
Z_k \ =\  \sum_{n\geq 1} I(L_n=k),\ \ \ \ {\rm for \ } k\geq 1.
\end{equation}
We note the zeroth mark $L_0$ is not included in the above
summation since any $Y_i$ with $i<L_0$ is part of an initial
segment of zeros of the sequence not preceded by a $1$, and so
does not contribute to any $d$-string, for $d\geq 1$.


\begin{theorem}
\label{cmpp_thm}
Suppose $\int \lambda_{w}(x) q(x,k) dx<\infty$ for all $w\in \R$
and $k\geq 1$.  Then, the count vector ${\bf Z}$ associated with
sequence ${\bf Y}$, defined through CMPP $({\bf X,L})={\mathcal
  M}(g,r,\lambda,q)$, is distributed as follows.  Given the value
$X_0=x_0$,
$${\bf Z}\ \stackrel{d}{=}\ \prod_{k\geq
    1}{\rm Po}\bigg(\int
\lambda_{x_0}(x) q(x,k) dx\bigg).$$
\end{theorem}

\begin{remark}
\label{cmpp_rmk}\rm
The distribution of ${\bf Z}$ does not depend on
the transition function $r$, consistent with the discussion of $L_0$
before the theorem.

Also, for a given $k\geq 1$, $Z_k$ is infinite with positive
probability exactly when there is a set $B$ such that $P(X_0\in
B)>0$ and $\int \lambda_{w}(x) q(x,k) dx=\infty$ for $w\in B$.
\end{remark}

{\it Proof of Theorem~\ref{cmpp_thm}.}
Recall the count vector representation (\ref{count_rep}).
Conditional on $X_0=x_0$, the point process $M$ on $\R\times \N$
given by $M(A\times \{k\}) = \sum_{i\geq 1} \delta_{X_i}(A)
I(L_i=k)$ is a Poisson process on $\R \times \N$ with intensity
function $\lambda_{x_0}(x) q(x,k)$ (cf. Proposition 4.10.1~(b) Resnick~(1994)).  Hence, it follows that, given
$X_0=x_0$, the variables $M(\R\times \{k\})= \sum_{n\geq 1}
I(L_n=k) = Z_k$ are independent Poisson variables with respective
means $\int \lambda_{x_0}(x)q(x,k) dx$, for $k\geq 1$.  \qed

\sect{The sequence ${\rm Bern}(a,b)$}
\label{holst_sect}
We now derive the count vector distribution for the sequence
${\rm Bern}(a,b)$ using a CMPP model.  Denote, as usual, for
$\alpha,\beta>0$, the Beta function
\begin{equation}
\label{beta}
B(\alpha,\beta) \ =\
\frac{\Gamma(\alpha)\Gamma(\beta)}{\Gamma(\alpha+\beta)},\end{equation}
and let
\begin{itemize}
\item [1.] $\bar{g}(x) = x^{b-1}(1-x)^{a-1}/B(b,a)$ on $0<x<1$, the
  Beta$(b,a)$ pdf,
\item [2.] $\bar{r}(x,k) = x^{k-1}(1-x)$ for $k\geq 1$,
\item [3.] $\bar{\lambda}_{w}(x) = [{a}/{(1-x)}]I(w<x<1)$, and
\item [4.] $\bar{q}(x,k) = x^{k-1}(1-x)$ for $k\geq 1$.
\end{itemize}

\begin{proposition}
\label{holst_prop}
The model $({\bf X,L})={\mathcal
  M}(\bar{g},\bar{r},\bar{\lambda},\bar{q})$ produces an
independent\break Bernoulli sequence ${\bf Y}\stackrel{d}{=}{\rm
  Bern}(a,b)$ for $a>0$ and $b>0$ whose count vector ${\bf Z}$, conditional on the
value $x_0$ of a ${\rm Beta}(b,a)$ random variable, is
distributed as $\prod_{k\geq 1}{\rm Po}(a(1-x_0^k)/k)$.

\end{proposition}

\begin{remark}
\label{holst_rmk}\rm
As a corollary, by taking $b \downarrow 0$, we recover the count
vector distribution for ${\rm Bern}(a,0)$ already considered in the
literature as simply ${\bf Z}\stackrel{d}{=}\prod_{k\geq 1} {\rm
  Po}(a/k)$.  Note that $(X_0,L_0) \rightarrow (0,1)$ in distribution
as $b \downarrow 0$.

The Poisson process in the above CMPP model with
intensity $\bar{\lambda}_{w}(\cdot)$ can be generated in the
following way.  First, the point process formed by the record values
from an iid sequence of
Beta$(1,a)$ random variables is a Poisson process with intensity
$a/(1-x)$, the Beta$(1,a)$ failure rate (cf. Resnick~(1994)
Proposition 4.11.1~(b)).  Next, we thin this process as follows.  Let
$X_0 \stackrel{d}{=} {\rm Beta}(b,a)$, and $\{X_i\}_{i\geq 1}$ be the
record values from an iid sequence of ${\rm Beta}(1,a)$ random
variables, subject to $X_i>X_0$ for $i\geq 1$.  Then, conditional on
$X_0=x_0$, the point process $\bar{N}$ defined by $\bar{N}(A) =
\sum_{i\geq 1}\delta_{X_i}(A)$ is the desired Poisson process with
intensity function $\bar{\lambda}_{x_0}(x)=[a/(1-x)]I(x_0<x<1)$.
\end{remark}

{\it Proof of Proposition~\ref{holst_prop}.}  The second part on the
count vector distribution follows from Theorem~\ref{cmpp_thm}, noting
for $k\geq 1$, that
\begin{equation}
\label{indep_mean_comp}
\int_0^1\bar{\lambda}_{x_0}(x)\bar{q}(x,k)dx \ = \ \int_{x_0}^1
ax^{k-1}dx \ = \ \frac{a(1-x_0^k)}{k}.\end{equation}

For the first part, we observe that the
distribution of $\{Y_i\}_{i\geq 1}$
given through (\ref{bernoulli_seq}) is uniquely determined by
the probabilities of cylinder sets of the form
\begin{eqnarray}
  \label{cylinder}
  &&E(k_0,\dots,k_n) = (L_0=k_0,L_1=k_1,\dots, L_n=k_n)\\\nonumber
  &=& \Big(Y_{t} = 1 {\rm \
    for \ } t \in
  \{K_0,K_1,\dots,K_n\}, {\rm and \ } Y_t = 0 {\rm \ otherwise \ for \ }
  1 \le t \le K_n\Big) \nonumber
 \end{eqnarray}
where $k_0,k_1,\dots,k_n$ are positive integers and
$K_0=k_0,K_1=K_0+k_1,\dots,K_n=K_{n-1}+k_n$ are their partial sums.
If the probability of
sets of the form $E\stackrel{def}{=}E(k_0,\dots,k_n)$ is a product of appropriate marginal
probabilities then $\{Y_n, n\ge 1\}$ will be the Bernoulli
sequence ${\rm Bern}(a,b)$. We will proceed to establish this.

Let $A_n =
\{0<x_0<x_1<\cdots<x_n<1\}$.  Using the Beta variables
representation
in Remark~\ref{holst_rmk}, write
\begin{eqnarray*}
P(E) &=&\int_{A_n} \bar{g}(x_0)\bar{r}(x_0,k_0) \prod_{i= 1}^n\Big[
P(X_i\in dx_i|X_i>x_{i-1}) \bar{q}(x_i,k_i)\Big]
dx_0.
\end{eqnarray*}
Since $P(X_i\in dx_i|X_i>x_{i-1}) =
a(1-x_i)^{a-1}/(1-x_{i-1})^a\,dx_i$ for $1\leq i\leq n$, we have
further that the last line equals
\begin{eqnarray}
\label{calc}
&&\frac{a^n}{B(b,a)} \int_{A_n} x_0^{b+k_0-2}\prod_{i= 1}^n x_i^{k_i-1}
(1-x_n)^a dx_0\dots dx_n \\
&&\ \ \ \ \ \ \ \ \ \ \ \ \ \ \ \ \ \ \ \ \ \ \ = \ \frac{B(b+K_n-1,a+1)}{B(b,a)}\cdot\frac{a^n}{\prod_{s=0}^{n-1}(b+K_s-1)}\nonumber
\end{eqnarray}
and, noting (\ref{beta}) and
$\alpha\Gamma(\alpha) = \Gamma(\alpha+1)$, that (\ref{calc}) becomes
\begin{eqnarray*}
\frac{a\prod_{r=0}^{K_n-2}(b+r)}{\prod_{r=0}^{K_n-1}(a+b+r)}\cdot
\frac{a^{n}}{\prod_{s=0}^{n-1}(b+K_s-1)}
& =& \ \prod_{i=1}^{K_n}\frac{b+i-1}{a+b+i-1} \prod_{r=0}^n \frac{a}{b+K_r
  -1}\end{eqnarray*}
which is exactly
$\prod_{i=1}^{K_n} P(Y_i =0) \prod_{r=0}^n
[P(Y_{K_r}=1)/P(Y_{K_r}=0)]$
with ${\bf Y}$ specified as ${\rm Bern}(a,b)$.
\qed

\sect{The sequence ${\rm Bern}_1(a,b)$}
\label{new_sect}
We will derive the count vector distribution for the sequence
${\rm Bern}_1(a,b)$, and show a dichotomy depending on whether
$b\geq 1$ or $b<1$.  We first consider the case where $a>0$ and
$b>1$.  Define
\begin{itemize}
\item [1.] $g^*(x) = x^{b-2}(1-x)^{a}/B(b-1,a+1)$ on $0<x<1$, \mbox{the
  Beta$(b-1,a+1)$ pdf,}
\item [2.] $r^*(x,1) = 1$,
\item [3.] $\lambda^*_{w}(x) = [{a}/{(1-x)}]I(w<x<1)$, and
\item [4.] $q^*(x,k) = x^{k-1}(1-x)$ for $k\geq 1$.
\end{itemize}

\begin{proposition}
\label{new_prop}
The CMPP model $({\bf X,L})={\mathcal M}(g^*,r^*,\lambda^*,q^*)$
produces an independent Bernoulli sequence ${\bf
  Y}\stackrel{d}{=}{\rm Bern}_1(a,b)$ for $a>0$ and $b>1$, and,
conditional on a Beta$(b-1,a+1)$ variable $X_0=x_0$, the
distribution of its count vector ${\bf Z}$ is $\prod_{k\geq
  1}{\rm Po}(a(1-x_0^k)/k)$.
\end{proposition}

\begin{remark}\label{new_rmk}
  \rm As a corollary, by taking $b\downarrow 1$, we find the count
  vector distribution for ${\rm Bern}_1(a,1)$ to be simply ${\bf
    Z}\stackrel{d}{=}\prod_{k\geq 1} {\rm Po}(a/k)$.  [In fact, ${\rm
    Bern}_1(a,1)$ coincides with the sequence ${\rm Bern}(a,0)$ mentioned earlier
in Remark \ref{holst_rmk}.]

  Also, we note the Poisson process in the above CMPP model with
  intensity $\lambda^*$ can be generated, as in
  Proposition~\ref{holst_prop}, by taking $X_0 \stackrel{d}{=} {\rm
    Beta}(b-1,a+1)$, and $\{X_i\}_{i\geq 1}$ as the sequence of
  records from an iid sequence of ${\rm Beta}(1,a)$ random variables,
  subject to the condition $X_1 > X_0$.
\end{remark}

{\it Proof of Proposition~\ref{new_prop}.}  We need only
establish the distribution of ${\bf Y}$, as the last statement
follows from Theorem~\ref{cmpp_thm} and the computation
(\ref{indep_mean_comp}).  The calculations are similar to the
proof of Proposition~\ref{holst_prop}.  Let
$k_0=1,k_1,k_2,\dots,k_n$ be positive integers, and
$K_0=k_0=1,K_1=K_0+k_1,\dots,K_n=K_{n-1}+k_n$ be their partial
sums. Recall the cylinder set defined in (\ref{cylinder}) and let
\begin{eqnarray*}
E_1 \ \stackrel{def}{=}\ E(1,k_1,\ldots,k_n) \ = \
(L_0=1,L_1=k_1,\dots, L_n=k_n),\end{eqnarray*}
and set $A_n =
\{0<x_0<x_1<\cdots <x_n<1\}$.
Write, using the construction in Remark~\ref{new_rmk}, that
\begin{eqnarray*}
P(E_1)
&=&\frac{1}{B(b-1,a+1)} \int_{A_n} \big[x_0^{b-2}
(1-x_0)^{a}\big]\cdot 1
\\
&&\ \ \ \ \ \ \  \times\prod_{i=1}^n \big[a(1-x_i)^{a-1}/(1-x_{i-1})^a\big]\big[x_i^{k_i-1}(1-x_i)\big]
dx_0\dots dx_n\\
&=&\frac{a^n}{B(b-1,a+1)} \int_{A_n} x_0^{b-2}\prod_{i=1}^n x_i^{k_i-1}
(1-x_n)^a dx_0\dots dx_n.
\end{eqnarray*}
Then, with (\ref{beta}) and $\alpha\Gamma(\alpha) = \Gamma(\alpha
+1)$, the last line equals
\begin{eqnarray*}
&&\frac{B(b+K_n-2, a+1)}{B(b-1,a+1)}\cdot
\frac{a^n}{(b-1)\prod_{s=1}^{n-1}(b+K_s-2)}\\
&&\ \ \ \ \ \ \ \ \ \ \ \ \ \ = \
\frac{\prod_{r=0}^{K_n-2}(b-1+r)}{\prod_{r=0}^{K_n-2}(a+b+r)}\cdot
 \frac{a^{n}}{(b-1)\prod_{s=1}^{n-1}(b+K_s-2)}\\
&&\ \ \ \ \ \ \ \ \ \ \ \ \ \ = \ \prod_{i=1}^{K_n-1}\frac{b+i-1}{a+b+i-1} \prod_{r=1}^n \frac{a}{b+K_r
  -2}\end{eqnarray*}
which is exactly
$P(Y_1=1)\prod_{i=2}^{K_n} P(Y_i =0) \prod_{r=1}^n [
P(Y_{K_r}=1)/P(Y_{K_r}=0)]$
with ${\bf Y}$ specified as ${\rm Bern}_1(a,b)$.
\qed

We now give the distribution of the count vector under ${\rm
  Bern}_1(a,b)$ for all $a>0$ and $b \ge 0$ by conditioning on the
location of the second $1$ in the sequence ${\bf Y}$.  Denote ${\bf
  Z}(a,b)$ as the count vector with respect to ${\rm Bern}_1(a,b)$ for
$a>0$ and $b\geq 0$.  Let ${\bf W_n}$ be the sequence whose $n$th
co-ordinate is $1$ and all the other co-ordinates are zero, for $n\geq
1$.  Let also
\[p_n \ = \ \left\{\begin{array}{rl}
\frac{a}{a+b}& \ {\rm for \ }n=2\\
\frac{a}{a+b+n-2}\prod_{r=0}^{n-3}\frac{b+r}{a+b+r}& \ {\rm for \ }
n\geq 3\end{array}\right.\]
be the probability that the second $1$ in ${\rm Bern}_1(a,b)$ occurs
at time $n\geq 2$, and note $\sum_{n\geq 2} p_n = 1$.
\begin{proposition}
\label{decomp_prop}
For $a>0$ and $b\geq 0$, we have
\begin{equation}
\label{general_b}
\mathcal{L}\left( {\bf Z}(a,b)\right) \ =\  \sum_{n\geq 2} p_n
\,\mathcal{L}\Big({\bf Z}(a,b+n-1) + {\bf W}_{n-1}\Big),
\end{equation}
and ${\bf
  Z}(a,b+n-1)$, conditional
on the value $x_0$ of a ${\rm Beta}(b+n-2,a+1)$ random variable, is
 distributed as $\prod_{k\geq 1}{\rm Po}(a(1-x_0^k)/k)$, for $b>0$ and
 $n\geq 2$.

\end{proposition}

\begin{remark}
  \rm The special case $b=0$ is interesting. The sequence ${\rm
    Bern}_1(a,0)$ is the independent sequence where $Y_1=Y_2=1$ and
  $P(Y_n=1) = a/(a+n-2)$ for $n\geq 3$.  That is, starting from time
  $n=2$, the sequence is ${\rm Bern}_1(a,1)={\rm Bern}(a,0)$.  Hence,
  by Proposition~\ref{holst_prop} (see Remark \ref{holst_rmk}), ${\bf Z}(a,0)$ is distributed as
  $\hat{{\bf Z}}+{\bf W}_1$ where $\hat{{\bf Z}}\stackrel{d}{=}
  \prod_{k\geq 1}{\rm Po}(a/k)$ is the count vector for ${\rm
    Bern}(a,0)$. This agrees with (\ref{general_b}), since $p_2=1$
  (when $b=0$) and ${\bf Z}(a,1)= \hat{\bf Z}$.
\end{remark}

{\it Proof of Proposition~\ref{decomp_prop}.} The distribution of
${\bf Z}(a,b)$ follows by conditioning on the first time that $Y_n =1$
for $n\ge 2$. The distributions of ${\bf Z}(a,b+n-1)$ are
completely specified by Proposition~\ref{new_prop} and
Remark~\ref{new_rmk}, since $b+n-1\geq 1$ for $n\geq 2$.  \qed

From (\ref{general_b}), it is not clear whether the
distribution of ${\bf Z}(a,b)$ is a mixture of product Poisson
factors or not for $0\le b < 1$. We show now that even the
first component $Z_1(a,b)$ is not a mixture of Poissons when
$0\le b < 1$.

\begin{proposition}
\label{cont_prop}
The distribution of $Z_1 \equiv Z_1(a,b)$, the count of\ $1$-strings in
the\break ${\rm Bern}_1(a,b)$ sequence, is not a mixture of Poissons
when $0\leq b<1$, that is, there is no measure $\mu$ on $[0,\infty)$
such that
\begin{equation}
  \label{mixture_cond}
E\Big[\exp\{tZ_1\}\Big] \ = \ \int_{[0,\infty)} e^{v(e^t-1)}d\mu(v).
\end{equation}
\end{proposition}

{\it Proof.}
It is well known that when (\ref{mixture_cond}) holds, the
variable $Z_1$ is over-dispersed, that is $O(Z_1) \stackrel{def}{=}
\textrm{Var}(Z_1) - E(Z_1) \ge 0$.  The proof now follows
by the expression for  $O(Z_1)$ in \rr{eq:overdisp} below.
Let ${\bf Y} = {\rm Bern}_1(a,b)$. Then,
\begin{equation}
  \label{decomp_Z}
  Z_1 \ =\  Y_2 +   \hat{Z}_1 = Y_2 + Y_2Y_3 + Z_1^{+}
\end{equation}
where $\hat{Z}_1 =\sum_{i\ge 2} Y_iY_{i+1}$ and $Z_1^{+} =
\sum_{i\ge 3} Y_i Y_{i+1}$, and the latter is independent of $Y_2$.
Furthermore $\hat{Z}_1$, $Z_1^{+}$ are the counts of strings of order
$1$ from ${\rm Bern}(a,b)$, ${\rm Bern}(a,b+1)$, respectively, and
their distributions are known from Proposition~\ref{holst_prop}. Hence, by easy calculations
\[E(\hat{Z_1}) = \frac{a^2}{(a+b)},\,  E(Z_1^{+}) =
\frac{a^2}{(a+b+1)},\,  E(\hat{Z_1}^2) =
\frac{a^3(a+1)}{(a+b)(a+b+1)} + \frac{a^2}{(a+b)}.\]
From the identities in (\ref{decomp_Z}), we  have
\[E(Z_1) = \frac{a(a+1)}{(a+b)},\; E(Z_1^2)=
\frac{a(a+1)}{(a+b)} + \frac{a^2(a+1)(a+2)}{(a+b)(a+b+1)}.\]
This leads to
\begin{equation}
  \label{eq:overdisp}
O(Z_1) \ =\  \frac{a^2(a+1)(b-1)}{(a+b)^2(a+b+1)}
\end{equation}
which is negative for  $b<1$, and positive for $b>1$. \qed

\sect{Some dependent Bernoulli sequences}
\label{dep_sect}
Two examples of dependent Bernoulli sequences, arising
in CMPP models with simple structures, whose count vector
distributions are mixtures of independent
Poisson factors are given.

\vskip .2cm
{\bf First Sequence.} For $a>0$ and $b> 0$, denote $P_{a,b}$ as the probability distribution of the CMPP
${\mathcal M}(\bar{g},\bar{r},\bar{\lambda},\bar{q})$ described in
Proposition~\ref{holst_prop} which gives rise to the Bernoulli
sequence ${\rm Bern}(a,b)$.
Let now $r^+(x,k) = kx^{k-1}(1-x)^2$ for $k\geq 1$.
Consider the associated CMPP model
${\mathcal M}(\bar{g},r^+,\bar{\lambda},\bar{q})$ with
$\bar{g},\bar{\lambda},\bar{q}$ the same as in
Proposition~\ref{holst_prop}. Denote the probability measure under this
model as $P^+=P^+_{a,b}$.

Note that $r^+(x,k) =k[\bar{r}(x,k) - \bar{r}(x,k+1)]$ where
$\bar{r}(x,k)=x^{k-1}(1-x)$.
Recall the cylinder set $E\stackrel{def}{=}E(k_0,\ldots, k_n)$ from
(\ref{cylinder})
where $k_0,k_1,\dots,k_n$ are positive
integers, and $K_0,K_1,\dots,K_n$ their partial sums. It is easy to see that
\begin{eqnarray*}
P^+(E) &=& k_0\Big[P_{a,b} \Big(E(k_0,\dots,k_n)\Big)
 - P_{a,b}\Big(E(k_0+1,k_1,\dots,k_n)\Big)\Big].
\end{eqnarray*}
%
From this expression, the distribution of ${\bf Y}$ can be recovered,
and shown to be not that of independent Bernoulli variables.  For
instance,
$$
  P^+(Y_1=1) \ =\  P_{a,b}(Y_1=1) - P_{a,b}(Y_1=0,Y_2=1)
\ =\
\frac{a(a+1)}{(a+b)(a+b+1)},
$$
and analogously
\[P^+(Y_2=1) \ = \
\frac{a^2(a+2) + 2ba(a+1)}{(a+b)(a+b+1)(a+b+2)}.
\]
Thus
\[
P^+(Y_1=1)P^+(Y_2=1) \  =\
\frac{a^2(a+1)(a^2+2a+2ba+2b)}{(a+b)^2(a+b+1)^2(a+b+2)},
\]
which does not match
\[P^+(Y_1=1,Y_2=1) \ = \ \frac{a^2 (a+2)}{(a+b)(a+b+1)(a+b+2)}\]
for $a,b>0$.

Finally, by Remark~\ref{cmpp_rmk},
we note the count vectors under $P_{a,b}$ and $P^+$ have the
same distribution, and by Proposition~\ref{holst_prop}
conditional
on the value of $x_0$ of a ${\rm Beta}(b,a)$ variable,
the count vectors
are distributed as $\prod_{k\geq 1}{\rm Po}(a(1-x^k_0)/k)$.
\vskip .2cm

{\bf Second Sequence.}  Consider $P_{1,0}$, the measure for the
CMPP model discussed in Example~\ref{ex_3} and
Remark~\ref{holst_rmk}, with respect to Bernoulli sequence ${\rm
  Bern}(1,0)$,
where $(X_0, L_0) \equiv (0,1)$,
$\{X_i\}_{i\geq 1}$ are the records from an iid
Uniform$[0,1]$ sequence, and $L_i$ are Geometric$(1-X_i)$ for
$i\geq 1$.

Let $P'$ stand for the measure under the ``switched'' CMPP model
where $(X_1,L_1)$ and $(X_2,L_2)$ are interchanged.  The
probabilities of ${\bf Y}$ on cylinder sets (cf. (\ref{cylinder}),
under $P'$, is given
by
\begin{eqnarray*}
 P'\Big(E(1,k_1,\ldots,k_n)\Big)
%
  &=&  P'(L_1=k_1, \ldots, L_n=k_n)\\
  &=&  P_{1,0}(L_2=k_1,
  L_1=k_2,\ {\rm and \ }L_i=k_i {\rm \ for \ }3\leq i\leq n)\end{eqnarray*}
for positive integers $k_0=1,k_1,\ldots, k_n$,
with $K_0=1,K_1=K_0+k_1,\ldots, K_n=K_{n-1} + k_n$
as their partial sums.  Under both models $P_{1,0}$ and $P'$, as
only two terms ($L_1,L_2$) exchange places, the associated count vectors are
the same, and by Proposition~\ref{holst_prop} distributed as
$\prod_{k\geq 1}{\rm Po}(1/k)$.

We now show that $\{Y_i\}_{i\geq 1}$ is not an independent
sequence under $P'$.
From the calculation in (\ref{calc}) with $(X_0,L_0)\equiv (0,1)$,
$Y_1\equiv 1$
and $\bar{r}(x,1)=1$ (take $b\downarrow 0$), and $a=1$, we can
write
\begin{eqnarray*}
P'(Y_2=1) &=& P_{1,0}(L_2=1) \ = \ \sum_{k\geq 1}
P_{1,0}(L_1=k,L_2=1)\\
&=& \sum_{k\geq 1}\int_{0<x_1<x_2<1} x_1^{k-1} (1-x_2) dx_1 dx_2 \ = \
1/4.
\end{eqnarray*}
Also,
\begin{eqnarray*}
P'(Y_2=1,Y_3=1)&=& P_{1,0}(L_1=1,L_2=1)\ = \  P_{1,0}(Y_2=1,Y_3=1) \ =
\ 1/6,\\
P'(Y_2=0,Y_3=1) &=& P_{1,0}(L_2=2)\\
&=& \sum_{k\geq 1} \int_{0<x_1<x_2<1} x_1^{k-1} x_2 (1-x_2) dx_1 dx_2
\ =\  5/36,
\end{eqnarray*}
which give $P'(Y_3=1)= 11/36$.
However, $P'(Y_2=1) P'(Y_3=1) = 11/144 \ \neq \ 1/6 = P'(Y_2=1,Y_3=1)$.

\bibliographystyle{plain}

\end{document}